\documentclass[12pt]{article}
\usepackage{bbm}
\usepackage{mathrsfs}
\usepackage{amsfonts}
\usepackage{amssymb,amsmath}
\def\cl{\centerline}

\def\vs{\vspace*}

\def\ni{\noindent}

\def\QED{\hfill$\Box$}
\numberwithin{equation}{section}
\newtheorem{theo}{Theorem}[section]
\newtheorem{defi}[theo]{Definition}

\newtheorem{lemm}[theo]{Lemma}

\newtheorem{corol}[theo]{Corollary}

\def\adddot{$\!\!\!${\bf.}\ \ }

\begin{document}
\cl {{\Large\bf Quantization on Generalized Heisenberg-Virasoro
Algebra}\footnote{Supported by NSF of China (No.11001046), the
Fundamental Research Funds for the Central Universities,
"Outstanding young teachers of Donghua University"
foundation.\\$^*$Corresponding author }}

 \vs{6pt}

\cl{Haibo Chen$^{\,1)}$,  Ran Shen$^{{\,2)}*}$, Jiangang Zhang$^{\,
3)}$}

\vs{6pt}

 \cl{\small
$^{1)}$,$^{2)}$College of Science, Donghua University, Shanghai,
201620, China}

\cl{\small Email:rshen@dhu.edu.cn}

\cl{\small $^{3)}$Department of Mathematics, Shanghai Normal
University, Shanghai, 200234, China}

\vs{6pt}

 {\small
\parskip .005 truein
\baselineskip 3pt \lineskip 3pt

\noindent{{\bf Abstract:} In a recent paper by the authors, Lie
bialgebra structures on generalized Heisenberg-Virasoro algebra
$\mathfrak{L}$ are considered. In this paper, the explicit formula
of the quantization on generalized Heisenberg-Virasoro algebra is
presented. \vs{6pt}

}\vs{6pt}

\noindent{\bf MR(2000) Subject Classification}\,\,\,\, 17B62, 17B05, 17B37, 17B66

\parskip .001 truein\baselineskip 8pt \lineskip 8pt

\vs{6pt}
\section{Introduction}
\setcounter{section}{1}\setcounter{equation}{0} In Hopf algebras or
quantum groups theory, there are two standard methods to yield new
bialgebras from old ones, one is twisting the product by a 2-cocycle
but keeping the coproduct unchanged, another is twisting the
coproduct by a Drinfel'd twist element but keeping the product
unchanged. Constructing quantization of Lie bialgebras is an
important method to produce new quantum groups
(cf.{\cite{D2},{\cite{D3}},{\cite{ES}}, etc). Drinfel'd in
{\cite{D5}} formulated a number of problems in quantum group theory,
including the existence of a quantization for Lie bialgebras. In the
paper {\cite{EK1}} Etingof and Kazhdan gave a positive answer to
some of Drinfel'd questions. In particular, they showed the
existence of quantizations for Lie bialgebras, namely, any classical
Yang-Baxter algebra can be quantized. Since then the interests in
quantizations of Lie bialgebras have been growing in the
mathematical literatures (e.g.,{\cite{EH,EK2,G,HW}}).

This Lie algebra is the universal central extension of the Lie
algebra of differential operators on a circle of order at most one,
which contains an infinite-dimensional Heisenberg subalgebra and the
Virasoro subalgebra. The natural action of the Virasoro subalgebra
on the Heisenberg subalgebra is twisted with a 2-cocycle. The
structure and representation theory for the generalized
Heisenberg-Virasoro algebra has been well developed (e.g.,{
\cite{ACKP, B, FO, JJ, SJ}}). The structure of the irreducible
highest weight modules and verma modules for the twisted
Heisenberg-Virasoro algebra are determined in {\cite{ACKP, B}}.

Recently, the Lie bialgebra structures on generalized
Heisenberg-Virasoro algebra $\mathfrak{L}$ was discussed in
{\cite{CSZ}}, which turned out the centerless  generalized
Heisenberg-Virasoro algebra $\overline{\mathfrak{L}}$ is triangular
coboundary. We note that  generalized Heisenberg-Virasoro algebra is
$\Gamma$ graded, where $\Gamma$ is an abelian group over a field
$\mathbb{F}$ of characteristic zero. For $\Gamma$ and $T$ is a
vector space over field $\mathbb{F}$, the {\it generalized
Heisenberg-Virasoro algebra}
$\mathfrak{L}:=\mathfrak{L}{(\Gamma)}$({\cite{LZ}}) is a Lie algebra
generated by $\{L_{x}=t^x\partial,I_{x}=t^x,C_L,C_I,C_{LI},x\in
\Gamma\}$, subject to the following relations:
\begin{eqnarray}\label{LB}\begin{array}{lll}
&&[L_{x},L_{y}] = (y-x) L_{x+y}+ \delta_{x+y,0} \frac{1}{12}
(x^3 - x)C_L,\\[4pt]&&
[I_{x},I_{y}] = y \delta_{x+y,0} C_I,\\[4pt]&&
[L_{x},I_{y}]= y I_{x+y} + \delta_{x+y,0} (x^2- x)C_L,\\[4pt]&&
[\mathfrak{L},C_L] = [\mathfrak{L},C_I] = [\mathfrak{L},C_{LI}] =
0.\end{array}
\end{eqnarray}
The Lie algebra $\mathfrak{L}$ has a generalized Heisenberg
subalgebra and a generalized Virasoro subalgebra interwined with a
2-cocycle. Set $\mathfrak{L}_x = \text{Span}_\mathbb{F}
\{L_{x},I_{x}\}$ for $x\in \Gamma \setminus\{0\}$, $\mathfrak{L}_0
=\text{Span}_\mathbb{F} \{L_{0},I_{0},C_L,C_I,C_{LI}\}$. Then
$\mathfrak{L}=\underset{x\in \Gamma}{\oplus} \mathfrak{L}_x $ is a
graded Lie algebra. Denote $C$ the center of $\mathfrak{L}$, then
$\mathcal{C}=\text{Span}_\mathbb{F}\{I_0,C_L,C_I,C_{LI}\}$. Denote
$\overline{\mathfrak{L}}=\mathfrak{L}/\mathcal {C}$, then
$\overline{\mathfrak{L}}$ is the centerless generalized
Heisenberg-Virasoro algebra.

The main result of this paper is the following theorem:

\begin{theo}\adddot
\label{main} We choose two distinguished elements
$h=\alpha^{-1}L_{0}$ and $e=I_{\alpha}$ with $\alpha\in \Gamma
\setminus \{0\}$, such that $[h,e]=e$ in $\overline{\mathfrak{L}}$,
there exists a structure of noncommutative and noncocommutative Hopf
algebra $(U(\overline{\mathfrak{L}})[[t]], m,\iota,\Delta,
S,\epsilon)$ on $U(\overline{\mathfrak{L}})[[t]]$ with
$U(\overline{\mathfrak{L}})[[t]]/tU(\overline{\mathfrak{L}})[[t]]\cong
U(\overline{\mathfrak{L}})$,
 which preserves the product and counit of
$U(\overline{\mathfrak{L}})[[t]]$ but with a comultiplication and
antipode defined by:

\begin{eqnarray*}
&&\Delta(L_\beta)=1\otimes L_\beta+ L_\beta\otimes (1-et)^{\alpha^{-1}\beta}+ \alpha h^{(1)}\otimes (1-et)^{-1}I_{\alpha+\beta}t,\\
&&\Delta(I_\alpha)= 1\otimes I_\alpha + I_\alpha \otimes (1-et),\\
&&S(L_\beta)=-(1-et)^{-\alpha^{-1}\beta} L_{\beta} +
(1-et)^{-\alpha^{-1}\beta}h^{[1]}_{-\alpha^{-1}\beta}I_{\alpha+\beta}t,\\
&&S(I_\alpha)=-(1-et)^{-1} I_\alpha.
\end{eqnarray*}.
\end{theo}
 \vskip10pt
\section{Preliminaries}
In this section, we summarize some basic definitions and results
concerning Lie bialgebra structures which will be used in the
following discussions. For a detailed discussion of this subject we
refer the reader to the literatures (see {\cite{CSZ}} and references
therein).

Let $\overline{\mathfrak{L}}$ be the centerless generalized
Heisenberg-Virasoro algebra and $U(\overline{\mathfrak{L}})$ the
universal enveloping algebra of $\overline{\mathfrak{L}}$. Then
$U(\overline{\mathfrak{L}})$ is equiped with a natural Hopf
algebraic structure
$(U(\overline{\mathfrak{L}}),m,\iota,\Delta_0,S_0,\epsilon)$, i.e.,
\begin{eqnarray}\label{2}
\Delta_{0}(X)=X\otimes 1+1\otimes X,& S_{0}(X)=-X, &\epsilon(X)=0,
\end{eqnarray}
where $\Delta_{0}$ is a comultiplication, $\epsilon$ is a counit and
$S_0$ is an antipode. In particular,
\begin{eqnarray}
\Delta_{0}(1)=1 \otimes 1 &\mathrm{and}& \epsilon(1)=S_0(1)=1.
\end{eqnarray}

The following result is due to W. Michaelis(see {\cite{M}}).
\begin{theo}\adddot
\label{main2}Let $L$ be a Lie algebra containing two linear
independent elements $a$ and $b$ satisfying $[a,b]=k b$ with $0\neq
k\in\mathbb{F}$. Set $$r=a\otimes b-b\otimes a$$ and define a linear
map
$$\Delta_r:L\rightarrow L\otimes L$$
by setting
$$\Delta_r(x)=x\cdot r=[x,a]\otimes b-b\otimes[x,a]+a\otimes[x,b]-[x,b]\otimes a,\ \ \ \forall \ x\in L.$$
Then $\Delta_r$ equips $L$ with structure of a triangular coboundary
Lie bialgebra.
\end{theo}
An algebra $L$ equipped with a classical Yang-Baxter $r$-matrix $r$
is called a {\it classical Yang-Baxter algebra}. It was shown in
{\cite{EK1}} that any classical Yang-Baxter algebra can be
quantized.
\begin{defi}\adddot \label{def2}\rm
Let $(H,m,\iota,\Delta_{0},S_{0},\epsilon)$ be a Hopf algebra over a
commutative ring $R$. A Drinfel'd twist $\mathcal {F}$ on $H$ is an
invertible element of $H\otimes H$ such that
\begin{eqnarray*}
&&(\mathcal {F}\otimes 1)(\Delta_{0}\otimes {\rm Id})(\mathcal
{F})=(1\otimes \mathcal {F})({\rm Id}\otimes\Delta_{0})(\mathcal
{F}),\\
 &&(\epsilon\otimes {\rm Id})(\mathcal {F})=1\otimes
1=({\rm Id}\otimes\epsilon)(\mathcal {F}).
\end{eqnarray*}
\end{defi}
The following result is well known (see e.g., [D1], [ES]).
\begin{lemm}\adddot\label{Legr}Let $(H,m,\iota,\Delta_{0},S_{0},\epsilon)$ be a Hopf algebra over a
commutative ring and  $\mathcal {F}$ be a Drinfel'd twist on $H$,
then $w=m({\rm Id}\otimes S_{0})(\mathcal {F})$ is invertible in $H$
with $w^{-1}=m(S_{0}\otimes {\rm Id})(\mathcal {F}^{-1})$. Moreover,
define $\Delta$ :$H\rightarrow H\otimes H$ and $S$ :$H\rightarrow H$
by
$$\begin{array}{llll}
\Delta(x)=\mathcal {F}\Delta_{0}(x)\mathcal {F}^{-1},&
S=wS_{0}(x)w^{-1},  \ \ \forall \ x\in H.
\!\!\!\!\!\!\!\!\!\!\!\!\!\!\!\!\!\!\!\!\!\!\!\!\end{array}$$ Then
 $(H,m,\iota,\Delta,S,\epsilon)$ is a new Hopf algebra, which is said to be the
 twisting of $H$ by the Drinfel'd twist $\mathcal {F}$.
\end{lemm}
Let $\mathbb{F}[[t]]$ be a ring of formal power series. Assume that
$L$ is a triangular Lie bialgebra with a classical Yang-Baxter
$r$-matrix $r$. Denote by $U(L)$ the universal enveloping algebra of
$L$, with the standard Hopf algebra structure $(U(
L),m,\iota,\Delta_{0},S_{0},\epsilon)$. Now consider the
topologically free $F[[t]]$-algebra $U(L)[[t]]$ (see [p.4]
{\cite{ES}} ), which can be viewed as an associative
$\mathbb{F}$-algebra of formal power series with coefficients in
$U(L)$. Naturally, $U(L)[[t]]$ is equiped with an induced Hopf
algebra structure arising from that on $U(L)$. By abuse of notation,
we denote it by $(U( L)[[t]],m,\iota,\Delta_{0},S_{0},\epsilon)$.

An algebra A equipped with a classical Yang-Baxter r-matrix r is
called a {\it classical Yang-Baxter algebra}. It is showed in
{\cite{EK1}} that any classical Yang- Baxter algebra can be
quantized.

For any element $x$ of a unital $R$-algebra ($R$ is a ring) and
$a\in R$, we set (see, e.g., {\cite{GZ}})
\begin{eqnarray*}
&&x^{(n)}_{a}:=(x+a)(x+a+1)\cdots(x+a+n-1),\\
&&x^{[n]}_{a}:=(x+a)(x+a-1)\cdots(x+a-n+1),
\end{eqnarray*}
and $x^{(n)}:=x^{(n)}_{0}$,  $x^{[n]}:=x^{[n]}_{0}$.

\begin{lemm}\adddot\label{Legr3}{\rm (see {\cite{G, GZ}})} For any element $x$ of a
unital $\mathbb{F}$-algebra, $a,b\in \mathbb{F}$, and  $r,s,t\in
\mathbb{Z}$, one has
\begin{eqnarray}\label{fir-e}
&&x^{(s+t)}_{a}=x^{(s)}_{a}x^{(t)}_{a+s}, \ \ \ \ \ x^{[s+t]}_{a}=
x^{[s]}_{a}x^{[t]}_{a-s},  \ \ \ \ \
x^{[s]}_{a}=x^{(s)}_{a-s+1},\label{fir-e}\\
&&\mbox{$\sum\limits_{s+t=r}$}\frac{(-1)^{t}}{s!t!}x^{[s]}_{a}x^{(t)}_{b}
={a-b\choose r}=
\frac{(a-b)\cdots(a-b-r+1)}{r!},\ \ \ \ \label{fir-e1}\\
&&\mbox{$\sum\limits_
{s+t=r}$}{\displaystyle\frac{(-1)^{t}}{s!t!}}x^{[s]}_{a}x^{[t]}_{b-s}
={\displaystyle{a-b+r-1\choose r}}=\frac{(a-b)\cdots(a-b+r-1)}{r!}.\
\ \ \ \ \ \ \ \label{fir-e2}
\end{eqnarray}
\end{lemm}

The following popular result will be frequently used in the third
part of this paper.
\begin{lemm}\adddot\label{Legr4}{\rm (see e.g., [Proposition 1.3(4)]{\cite{SF}})} For any elements $x, y$ of an
associative algebra $A$, and $m\in\mathbb{Z}_{+}$, one has
\begin{eqnarray}
xy^{m}=\mbox{$\sum\limits_{k=0}^{m}$}(-1)^k{m\choose k}y^{m-k}({\rm
ad\,}y)^k(x).
\end{eqnarray}
\end{lemm}
\section{Proof of main results}
In this section, assume that $\mathfrak{L}$ is the generalized
Heisenberg-Virasoro algebra defined in (\ref{LB}). We have obtained
that the Lie bialgebra structures on centerless generalized
Heisenberg-Virasoro algebra $\overline{\mathfrak{L}}$ are triangular
coboundary, namely, there always exist solutions of CYBE in
$\overline{\mathfrak{L}}\otimes\overline{\mathfrak{L}}$. Therefore,
$\overline{\mathfrak{L}}$ can be quantized by the above arguments.
In what follows we will use a Drinfel'd twist (see Definition
\ref{def2}) to proceed the quantization on centerless generalized
Heisenberg-Virasoro algebra $\overline{\mathfrak{L}}$.

To describe quantizations of $U(\overline{\mathfrak{L}})$, we need
to construct explicitly Drinfel¡¯d twists according to Lemma
\ref{Legr}. Set
$$h:=\alpha^{-1}L_0,\ \ \ \ \ \ e:=I_\alpha$$ for a  fixed $\alpha\in\Gamma \backslash \{0\}$. It
is easily to see $[h,e]=e$ by (\ref{LB}). Then it follows from
Theorem \ref{main2} that  $r=h\otimes e-e\otimes h$ is a solution of
CYBE, namely, $r$ is a classical $r$-matrix. Now we can use this
$r$-matrix determined by $e$ and $h$ to construct a Drinfel'd twist.
This will be done by several lemmas.
\begin{lemm}\adddot\label{lemm1}
For $a\in\mathbb{F}$, $i\in\mathbb{Z}_{+}$, $n\in\mathbb{Z}$,
$\beta\in \Gamma$ and $\alpha\in \Gamma\backslash \{0\}$, one has
\begin{eqnarray*}\begin{array}{lll}
&&L_\beta h^{(i)}_a=h^{(i)}_{a-\alpha^{-1}\beta}L_\beta,  \ \ \ \ \
L_\beta h^{[i]}_{a}=h^{[i]}_{a-\alpha^{-1}\beta}L_\beta,
\\
&&I_\alpha h^{(i)}_a=h^{(i)}_{a-1}I_\alpha,     \ \ \ \ \
 I_\alpha h^{[i]}_a=h^{[i]}_{a-1}I_\alpha, \\
 &&e^nh^{(i)}_{a}=h^{(i)}_{a-n}e^n,\ \ \ \ \ \ \
e^nh^{[i]}_{a}=h^{[i]}_{a-n}e^n.
\end{array}
\end{eqnarray*}
\end{lemm}
\ni{\it Proof.} We only prove the first equation (the others can be
obtained similarly). Since $L_\beta h-h L_\beta =-\alpha^{-1}\beta
L_\beta$, there is nothing to prove for $i=1$. For the induction
step, suppose that it holds for $i$, then one has
\begin{eqnarray*}
L_\beta h^{(i+1)}_a&=&L_\beta
h^{(i)}_a(h+a+i)\\&=&h^{(i)}_{a-\alpha^{-1}\beta} L_\beta
(h+a+i)\\&=&h^{(i)}_{a-\alpha^{-1}\beta}(h-\beta+a+i)L_\beta
\\&=&h^{(i+1)}_{a-\alpha^{-1}\beta}L_\beta.
\end{eqnarray*}

Now for $a\in\mathbb{F}$, set
\begin{eqnarray}\label{fom1}
&\mathcal
{F}_a=\sum\limits^{\infty}_{r=0}\frac{(-1)^r}{r!}h^{[r]}_a\otimes
e^r t^r,& F_a
=\mbox{$\sum\limits^{\infty}_{r=0}$}\frac{1}{r!}h^{(r)}_a\otimes e^r
t^r,
\\
&u_a=m\cdot (S_0\otimes{ \rm Id})(F_a), &v_a=m\cdot({\rm Id}\otimes
S_0)(\mathcal {F}_a).\nonumber
\end{eqnarray}
Write $\mathcal {F}=\mathcal {F}_0$, $F=F_0$, $u=u_0$, $v=v_0$.
Since $S_0(h^{(r)}_a)=(-1)^rh^{[r]}_{-a}$ and $S_0(e^r)=(-1)^r e^r$,
we have
\begin{eqnarray}\label{fom2}
u_a=\mbox{$\sum\limits^{\infty}_{r=0}$}\frac{(-1)^r}{r!}h^{[r]}_{-a}e^r
t^r,&&v_a= \mbox{$\sum\limits^{\infty}_{r=0}$}\frac{1}{r!}h^{[r]}_a
e^r t^r.
\end{eqnarray}
\begin{lemm}\adddot\label{lemm2}
For $a,b\in\mathbb{F}$, one has
\begin{eqnarray*}
\mathcal {F}_aF_b=1\otimes(1-et)^{a-b},&& v_au_b=(1-et)^{-(a+b)}.
\end{eqnarray*}
\end{lemm}
\ni{\it Proof.} Using the equations (\ref{fir-e1}) and (\ref{fom1}),
we have
\begin{eqnarray*}
\mathcal {F}_aF_b&=&
\mbox{$\sum\limits^{\infty}_{r,s=0}$}\frac{(-1)^r}{r!s!}h^{[r]}_ah^{(s)}_b\otimes
e^r e^s t^r t^s\\&=&\mbox{$\sum\limits^{\infty}_{m=0}$}(-1)^m
\big(\mbox{$\sum\limits_{r+s=m}$}\frac{(-1)^s}{r!s!}h^{[r]}_{a}h^{(s)}_b\big)\otimes e^m t^m\\
&=&\mbox{$\sum\limits^{\infty}_{m=0}$}(-1)^m{a-b \choose m} \otimes
e^m t^m\\&=&1\otimes(1-et)^{a-b}.
\end{eqnarray*}
From (\ref{fir-e2}), (\ref{fom2}) and Lemma \ref{lemm1}, we obtain
that
\\[4pt]
\hspace*{60pt}$v_au_b=\sum\limits^{\infty}_{r,s=0}{\displaystyle\frac{(-1)^s}{r!s!}}h^{[r]}_a
e^r h^{[s]}_{-b} e^s
t^{r+s}$\\[4pt]
\hspace*{60pt}$\phantom{v_au_b}=\sum\limits^{\infty}_{m=0}\sum\limits_{r+s=m}{\displaystyle\frac{(-1)^s}{r!s!}}
h^{[r]}_{a}h^{[s]}_{-b-r}e^m t^m$
\\[4pt]
\hspace*{60pt}$\phantom{v_au_b}=\sum\limits^{\infty}_{m=0}{\displaystyle{a+b+m-1\choose
m}}e^m t^m$\\[4pt]
\hspace*{60pt}$\phantom{v_au_b}=(1-et)^{-(a+b)}.$
\begin{corol}\adddot\label{coro}  For any $a\in\mathbb{F}$, the elements $F_a$ and $u_a$ are
invertible with $F^{-1}_a=\mathcal {F}_a$, $u^{-1}_a=v_{-a}$. In
particular, $F^{-1}=\mathcal{F}$, $u^{-1}=v$.
\end{corol}
\begin{lemm}\adddot\label{lemm3}
For any $a\in\mathbb{F}$ and $r\in\mathbb{Z}_+$, one has $
\Delta_0(h^{[r]})=\mbox{$\sum\limits^r_{i=0}$}{r\choose
i}h^{[i]}_{-a}\otimes h^{[r-i]}_{a}. $
 In particular, one has $\Delta_0(h^{[r]})=\mbox{$\sum\limits^r_{i=0}$}{r\choose i}h^{[i]}\otimes
 h^{[r-i]}$.
\end{lemm}
\ni{\it Proof.} Since $\Delta_0(h)=1\otimes h+h\otimes 1$, it is
easy to see that the result is true for $r=1$. Suppose it is true
for $r$, then for $r+1$, we have
\begin{eqnarray*}
\Delta_0(h^{[r+1]})&=&\Delta_0(h^{[r]}(h-r))\\
&=&\Delta_0(h^{[r]})\big(\Delta_0(h)-\Delta_0(r)\big)\\
&=&\Big(\mbox{$\sum\limits^r_{i=0}$}{r\choose i}h^{[i]}_{-a}\otimes
h^{[r-i]}_a\Big)\big((h-r)\otimes 1+1\otimes (h-r)+r(1\otimes 1)\big)\\
&=&\Big(\mbox{$\sum\limits^{r-1}_{i=1}$}{r\choose
i}h^{[i]}_{-a}\otimes h^{[r-i]}_a\Big)\big((h-r)\otimes 1+1\otimes
(h-r)\big)+r\mbox{$\sum\limits^r_{i=0}$}{r\choose
i}h^{[i]}_{-a}\otimes h^{[r-i]}_a \\&& +h^{[r+1]}_{-a}\otimes
1+h^{[r]}_{-a}\otimes a+h^{[r]}_{-a}\otimes(h-r)+(h-r)\otimes
h^{[r]}_a+1\otimes
h^{[r+1]}_a-a\otimes h^{[r]}_a\\
&=&1\otimes h^{[r+1]}_a+h^{[r+1]}_{-a}\otimes
1+r\mbox{$\sum\limits^{r-1}_{i=1}$}{r\choose i}h^{[i]}_{-a}\otimes
h^{[r-i]}_a+h^{[r]}_{-a}\otimes (h+a)\\&&
 +(h-a)\otimes h^{[r]}_a
+\mbox{$\sum\limits_{i=1}^{r-1}$}{r\choose i}h^{[i+1]}_{-a}\otimes
h^{[r-i]}_a+\mbox{$\sum\limits_{i=1}^{r-1}$}(-r+a+i){r\choose
i}h^{[i]}_{-a}\otimes
h^{[r-i]}_a\\
&&+\mbox{$\sum\limits_{i=1}^{r-1}$}{r\choose i}h^{[i]}_{-a}\otimes
h^{[r-i+1]}_a+\mbox{$\sum\limits_{i=1}^{r-1}$}(-a-i){r\choose
i}h^{[i]}_{-a}\otimes h^{[r-i]}_a\\
&=&1\otimes h^{[r+1]}_{a}+h^{[r+1]}_{-a}\otimes
1+\mbox{$\sum\limits_{i=1}^{r}$}\bigg[{r\choose i-1}+{r\choose
i}\bigg]h^{[i]}_{-a}\otimes h^{[r-i+1]}_a\\
&=&\mbox{$\sum\limits_{i=0}^{r+1}$}{r+1 \choose
i}h^{[i]}_{-a}\otimes h^{[r+1-i]}_a.
\end{eqnarray*}
Hence, the formula holds by induction.

\begin{lemm} \adddot \label{lemm5}
The element $\mathcal
{F}=\mbox{$\sum\limits^{\infty}_{r=0}$}\frac{(-1)^r}{r!}h^{[r]}\otimes
e^rt^r$ is a Drinfeld twist on $U(\overline{\mathfrak{L}})[[t]]$ .
\end{lemm}
\ni{\it Proof.} It can be proved directly by the similar methods as
those presented in the proof of [Proposition 2.5]{\cite{HW}}.

Now we can now perform the process of twisting the standard Hopf
structure $(U(\overline{\mathfrak{L}}),m,\iota,\Delta_0,
S_0,\epsilon)$ defined in (\ref{2}) by the Drinfel'd twist $\mathcal
{F}$ constructed above. The following lemmas are very useful to our
main results.
\begin{lemm}\adddot \label{lemm6} For $a\in\mathbb{F}$, $\beta\in \Gamma$, $\alpha\in \Gamma \backslash \{0\}$,
one has
\begin{eqnarray*}\begin{array}{lll}
(L_\beta\otimes 1)F_a=F_{a-\alpha^{-1}\beta}(L_\beta \otimes 1),\\
(I_\alpha \otimes 1)F_a=F_{a-1}(I_\alpha \otimes 1).
\end{array}
\end{eqnarray*}
\end{lemm}
\ni{\it Proof.} It follows directly from equation (\ref{fom1}) and
Lemma \ref{lemm1}.
\begin{lemm}\adddot\label{lemm7}
For $a\in\mathbb{F}$, $\beta\in \Gamma$, $\alpha\in \Gamma\backslash
\{0\}$ and $r\in\mathbb{Z}_+$, one has
\begin{eqnarray}
&&L_\beta e^r=e^rL_\beta +\alpha re^{r-1}I_{\alpha+\beta},\ \ \ \ \ \ \ \label{fom4}\\[6pt]
&&I_\alpha e^r=e^r I_\alpha.\ \ \ \ \ \ \ \ \label{fom5}
\end{eqnarray}
\end{lemm}
\ni{\it Proof.} By Lemma \ref{Legr4} and equation (\ref{LB}), we
have
\begin{eqnarray*}
L_\beta e^r&=&\mbox{$\sum\limits^r_{i=0}$}(-1)^i{r\choose
i}e^{r-i}({\rm
ad\,} e)^{i}(L_\beta)\\
&=&e^rL_\beta +\alpha re^{r-1}I_{\alpha+\beta}.
\end{eqnarray*}
 Similarly, one can get (\ref{fom5}).
\begin{lemm}\adddot\label{lemm8}  For $a\in\mathbb{F}$, $\beta\in \Gamma,\alpha\in
\Gamma \backslash \{0\}$, we have
\begin{eqnarray}
(1\otimes L_\beta)F_a &=&F_a(1\otimes L_\beta)+\alpha
F_{a+1}(h_a^{(1)}\otimes
I_{\alpha + \beta}t),\label{fom7}\\
(1\otimes I_\alpha)F_a &=&F_a(1\otimes I_\alpha).\label{fom9}\
\end{eqnarray}
\end{lemm}
\ni{\it Proof.} By equations (\ref{fir-e}), (\ref{fom1}) and
(\ref{fom4}), one has
\\[4pt]\hspace*{4ex}$
 (1\otimes L_\beta)F_a=\mbox{$\sum\limits^{\infty}_{r=0}$}\frac{1}{r!}h^{(r)}_a\otimes L_\beta e^r
 t^r$\\[4pt]\hspace*{4ex}$\phantom{ (1\otimes L_n)F_a}
 =\mbox{$\sum\limits^{\infty}_{r=0}$}\frac{1}{r!}h_a^{(r)}\otimes (e^rL_\beta +\alpha re^{r-1}I_{\alpha+\beta})t^r$\\[4pt]
\hspace*{4ex}$\phantom{ (1\otimes L_\beta)F_a}
 =F_a(1\otimes
 L_\beta)+\mbox{$\sum\limits^{\infty}_{r=1}$}\frac{\alpha}{(r-1)!}
h^{(r)}_{a}\otimes e^{r-1}I_{\alpha+\beta}t^r
$\\[4pt]\hspace*{4ex}$\phantom{ (1\otimes L_\beta)F_a}
 =F_a(1\otimes
 L_\beta)+\mbox{$\sum\limits^{\infty}_{r=0}$}\frac{\alpha}{r!}
h^{(r+1)}_{a}\otimes
e^{r}I_{\alpha +\beta}t^{r+1}$\\[4pt]\hspace*{4ex}$\phantom{ (1\otimes L_\beta)F_a}
 =F_a(1\otimes
 L_n)+ \mbox{$\sum\limits^{\infty}_{r=0}$}\frac{\alpha}{r!}
h^{(r)}_{a+1}h_a^{(1)}\otimes
e^{r}I_{\alpha+\beta}t^{r+1} $\\[4pt]\hspace*{4ex}$\phantom{ (1\otimes L_\beta)F_a}
 =F_a(1\otimes L_n)+ \alpha F_{a+1}\big(h^{(1)}_a\otimes
 I_{\alpha+\beta} t)$.\\
This proves equation (\ref{fom7}). Similarly,
 (\ref{fom9}) follow from
(\ref{fom5}) .
\begin{lemm}\adddot\label{lemm9}
 For $a\in\mathbb{F}$, $\beta\in \Gamma,\alpha\in \Gamma \backslash
\{0\}$,
 one has
 \begin{eqnarray}
 &&L_\beta u_a=u_{a+\alpha^{-1}\beta}L_\beta- u_{a+\alpha^{-1}\beta} h^{[1]}_{-a-\alpha^{-1}\beta}I_{\alpha+\beta} t,\label{fom10}\\
 &&I_\alpha u_a=u_{a+1}I_\alpha.\label{fom12}
 \end{eqnarray}
\end{lemm}
\ni{\it Proof.} From equations (\ref{fir-e}), (\ref{fom2}),
(\ref{fom4}) and Lemma \ref{lemm1}, one has
\\[4pt]\hspace*{4ex}$
 L_\beta u_a=\mbox{$\sum\limits^{\infty}_{r=0}$}\frac{(-1)^r}{r!}L_\beta h^{[r]}_{-a}e^rt^r$\\
 [4pt]\hspace*{4ex}$\phantom{L_\beta u_a}
 =\mbox{$\sum\limits^{\infty}_{r=0}$}
 \frac{(-1)^r}{r!}h^{[r]}_{-a-\alpha^{-1}\beta}L_\beta e^rt^r$\\
 [4pt]\hspace*{4ex}$\phantom{L_\beta u_a}
 =\mbox{$\sum\limits^{\infty}_{r=0}$}\frac{(-1)^r}{r!}h^{[r]}_{-a-\alpha^{-1}\beta}\big(e^rL_\beta+\alpha re^{r-1}
 I_{\alpha+\beta})t^r$\\
[4pt]\hspace*{4ex}$\phantom{L_\beta u_a}
=u_{a+\alpha^{-1}\beta}L_\beta+\mbox{$\sum\limits^{\infty}_{r=1}$}\frac{(-1)^r\alpha}{(r-1)!}h^{[r]}_{-a-\alpha^{-1}\beta}e^{r-1}
I_{\alpha+\beta}t^{r}$\\
[4pt]\hspace*{4ex}$\phantom{L_\beta u_a}
=u_{a+\alpha^{-1}\beta}L_\beta-\mbox{$\sum\limits^{\infty}_{r=0}$}\frac{(-1)^r}{r!}h^{[r+1]}_{-a-\alpha^{-1}\beta}e^{r}I_{\alpha+\beta}t^{r+1}$\\
[4pt]\hspace*{4ex}$\phantom{L_\beta u_a}
=u_{a+\alpha^{-1}\beta}L_\beta-\mbox{$\sum\limits^{\infty}_{r=0}$}
\frac{(-1)^r}{r!}h^{[r]}_{-a-\alpha^{-1}\beta}h^{[1]}_{-a-\alpha^{-1}\beta-r}e^{r} I_{\alpha+\beta}t^{r+1}$\\
[4pt] \hspace*{4ex}$\phantom{L_\beta u_a}
=u_{a+\alpha^{-1}\beta}L_\beta-
\mbox{$\sum\limits^{\infty}_{r=0}$}\frac{(-1)^r}{r!}h^{[r]}_{-a-\alpha^{-1}\beta}e^r
h^{[1]}_{-a-\alpha^{-1}\beta} I_{\alpha+\beta} t^{r+1}$\\
[4pt]\hspace*{4ex}$\phantom{L_\beta u_a}
=u_{a+\alpha^{-1}\beta}L_\beta-u_{a+\alpha^{-1}\beta}h^{[1]}_{-a-\alpha^{-1}\beta}I_{\alpha+\beta}t.$\\
 Hence, (\ref{fom10}) holds. Similarly,
one can get (\ref{fom12})  by Lemma \ref{lemm1} and
 Lemma \ref{lemm7}.

Now we have enough in hand to prove our main theorem in this paper.

{\it Proof of Theorem \ref{main}.}\ \ By Lemma
\ref{Legr}, Lemma\ref{lemm2}, Corollary \ref{coro}, Lemma\ref{lemm6} and Lemma\ref{lemm8}, we have\\
[4pt]\hspace*{4ex}$
 \Delta(L_\beta)=\mathcal
{F}\Delta_0(L_n)\mathcal{F}^{-1}$\\
[4pt]\hspace*{4ex}$\phantom{\Delta(L_\beta)}
 =\mathcal
{F}(L_\beta\otimes
1)F+\mathcal {F}(1\otimes L_\beta)F$\\
[4pt]\hspace*{4ex}$\phantom{\Delta(L_\beta)}
=\mathcal
{F}F_{-\alpha^{-1}\beta} (L_\beta \otimes
1)+\mathcal{F}\big(F(1\otimes L_\beta)+\alpha F_1
(h^{(1)}\otimes I_{\alpha+\beta})t)$\\
[4pt]\hspace*{4ex}$\phantom{\Delta(L_\beta)}
 =L_\beta\otimes (1-et)^{\alpha^{-1}\beta}+1\otimes
L_\beta+ \alpha h^{(1)}\otimes (1-et)^{-1}I_{\alpha+\beta}t
$\\[4pt]\hspace*{4ex}$\phantom{\Delta(L_\beta)}
=1\otimes L_\beta+ L_\beta\otimes (1-et)^{\alpha^{-1}\beta}+\alpha
h^{(1)}\otimes (1-et)^{-1}I_{\alpha+\beta}t
.$\\[4pt]\hspace*{4ex}$
\Delta(I_\alpha)=\mathcal {F}\Delta_0(I_\alpha)\mathcal{F}^{-1}$\\
[4pt]\hspace*{4ex}$\phantom{\Delta(I_\alpha)}
=\mathcal
{F}(I_\alpha\otimes 1)F+\mathcal {F}(1\otimes I_\alpha)F$\\
[4pt]\hspace*{4ex}$\phantom{\Delta(I_\alpha)} =\mathcal
{F}F_{-1}(I_\alpha\otimes 1)+\mathcal{F} F(1\otimes
I_\alpha)$\\
[4pt]\hspace*{4ex}$\phantom{\Delta(I_\alpha)}
 =1\otimes I_\alpha+
I_\alpha\otimes (1-et).
$\\[4pt]\hspace*{4ex}$\phantom{\Delta(I_\alpha)}$\\
Again by Lemma \ref{Legr}, Lemma\ref{lemm2},  Corollary \ref{coro}
and Lemma \ref{lemm9}, we have
\begin{eqnarray*}
S(L_\beta)&=&u^{-1}S_0(L_\beta)u\\
&=&-vL_\beta u\\
&=&-v(u_{\alpha^{-1}\beta} L_\beta- u_{\alpha^{-1}\beta}
h^{[1]}_{-\alpha^{-1}\beta}I_{\alpha+\beta}t)\\
&=&-(1-et)^{-\alpha^{-1}\beta}L_\beta+(1-et)^{-\alpha^{-1}\beta} h^{[1]}_{-\alpha^{-1}\beta}I_{\alpha+\beta}t.\\
S(I_\alpha)&=&u^{-1}S_0(I_\alpha)u\\
&=&-vI_\alpha u\\
&=&-v u_1I_\alpha\\
&=&-(1-et)^{-1}I_\alpha.\\
\end{eqnarray*}
So the proof is complete!

.\QED

\vskip12pt

%Included for Gather Purpose only:
%input "Xbib.bib"
\bibliographystyle{amsplain}

\end{document}